\documentclass[12pt]{amsart}
\usepackage[english]{babel}
\usepackage[latin1]{inputenc}
\usepackage{amstext}
\usepackage{amsmath}
\usepackage{amsfonts}
\usepackage{amssymb}
\usepackage{amsthm}
\usepackage[all]{xy}
\xyoption{arc}
\CompileMatrices
\usepackage{graphicx}

\newtheorem{thm}{Theorem}[section]

\newtheorem{lem}[thm]{Lemma}
\newtheorem{prop}[thm]{Proposition}
\newtheorem{dfn}[thm]{Definition}
\newtheorem{rem}[thm]{Remark}

\newenvironment{pr}[1][Proof]{\noindent\textbf{#1.} }{\ \rule{0.5em}{0.5em}}

\newcommand{\Hom}{\operatorname{Hom}}

\newcommand{\Proj}{\operatorname{Proj}}

\selectfont

\numberwithin{equation}{section}

\begin{document}

\title[Harder-Narasimhan filtration for quivers] {On the Harder-Narasimhan filtration for finite dimensional representations of quivers}
\author[A. Zamora]{Alfonso Zamora}

\address{Instituto de Ciencias Matem\'aticas (CSIC-UAM-UC3M-UCM),
Nicol\'as Cabrera 13-15, Campus Cantoblanco UAM, 28049 Madrid,
Spain}

\address{Departamento de \'{A}lgebra, Facultad de Matem\'aticas, Universidad Complutense de
Madrid, 28040 Madrid, Spain}

\email{alfonsozamora@icmat.es, alfonsozamora@mat.ucm.es}

\begin{abstract}
We prove that the Harder-Narasimhan filtration for an unstable
finite dimensional representation of a finite quiver coincides
with the filtration associated to the 1-parameter subgroup of
Kempf, which gives maximal unstability in the sense of Geometric
Invariant Theory for the corresponding point in the parameter
space where these objects are parametrized in the construction of
the moduli space.

$$ $$


\end{abstract}

\maketitle

\section*{Introduction}
Let $Q$ be a finite quiver, given by a finite set of vertices and
arrows between them, and a representation of $Q$ on finite
dimensional $k$-vector spaces, where $k$ is an algebraically
closed field of arbitrary characteristic. There exists a notion of
stability for such representations given by King (\cite{Ki}) and,
more generally by Reineke (\cite{Re}) (both particular cases of
the abstract notion of stability for an abelian category that we
can find in \cite{Ru}), and a notion of the existence of a unique
Harder-Narasimhan filtration with respect to that stability
condition.

We consider the construction of a moduli space for these objects
by King (\cite{Ki}) and associate to an unstable representation an
unstable point, in the sense of Geometric Invariant Theory, in a
parameter space where a group acts. Then, the $1$-parameter subgroup given by Kempf
(\cite{Ke}), which is maximally destabilizing in the GIT sense,
gives a filtration of subrepresentations and we prove that it
coincides with the Harder-Narasimhan filtration for that
representation.

This article makes use of the same techniques that a previous work of the
author in collaboration with T. G\'omez and I. Sols (\cite{GSZ}).
In that article, we considered an unstable torsion free sheaf $E$
over a smooth projective variety $X$. There, we proved that the
filtration associated to the $1$-parameter subgroup given by
Kempf, coincides with the Harder-Narasimhan filtration of $E$ with
the definition of stability given by Gieseker.

The definition of stability for a representation of a quiver (c.f.
Definition \ref{Qstability}) contains two sets of parameters,
the coefficients of the linear functions $\Theta$ and $\sigma$. In
\cite{Ke}, the $1$-parameter subgroup is taken to maximize certain
function which depends on the choice of a linearization of the action of 
the group we are taking the quotient by, 
and a \emph{length} in the set of $1$-parameter subgroups
(c.f. Definition \ref{length}). In the case of sheaves the group
is $SL(N)$, which is simple, so any such length is unique
up to multiplication by a scalar, whereas for finite dimensional
representations of quivers we quotient by a product of general
linear groups, so we have to choose a scalar for each factor in
the choice of a length. Hence, we put the positive coefficients of $\sigma$ precisely as these scalars
and consider a particular linearization depending on $\sigma$ and $\Theta$,
in order to relate the Harder-Narasimhan filtration of a representation with the filtration given in \cite{Ke} (c.f. Theorem \ref{KempfHN}).

Hesselink shows in \cite{He} that the unstable locus of a smooth complex projective variety acted by a complex reductive group can be stratified by conjugacy classes 
of $1$-parameter subgroups. Tur shows in \cite{Tu} that Hesselink's stratification coincides with a stratification in Harder-Narasimhan types (meaning the numerical invariants
appearing on the Harder-Narasimhan filtration), for the space of quiver representations. In \cite{Ho}, Hoskins proves that these two stratifications do coincide 
with a Morse stratification given by the norm square of the moment map, for reductive groups acting on affine spaces, in particular quiver representations.  

{\bf Acknowledgments.} The author wishes to thank L. \'{A}lvarez-C\'onsul, T. G\'omez and A. D. King, for useful discussions and comments. 
This work has been supported by project MTM2010-17389 and ICMAT Severo Ochoa project SEV-2011-0087
granted by Spanish Ministerio de Econom\'ia y Competitividad. The author was also supported by a FPU grant from the Spanish Ministerio de Educaci\'on.

\section{Harder-Narasimhan filtration for representations of quivers}
A finite quiver $Q$ is given by a finite set of vertices
$Q_{0}$ and a finite set of arrows $Q_{1}$. The arrows will be
denoted by $(\alpha: v_{i}\rightarrow v_{j})\in Q_{1}$. We denote
by $\mathbb{Z}Q_{0}$ the free abelian group generated by $Q_{0}$.

Fix $k$, an algebraically closed field of arbitrary
characteristic. Let $\mod kQ$ be the category of
finite dimensional representations of $Q$ over $k$. Such category
is an abelian category and its objects are given by tuples
$$M=((M_{v})_{v\in Q_{0}},(M_{\alpha}:M_{v_{i}}\rightarrow M_{v_{j}})_{\alpha:v_{i}\rightarrow v_{j}})$$
of finite dimensional $k$-vector spaces and $k$-linear maps
between them. The dimension vector of a representation is given by
$\underline{\dim}M=\sum_{v\in Q_{0}}\dim_{k}M_{v}\cdot
v\in\mathbb{N}Q_{0}$.

Let $\Theta$ be a set of numbers $\Theta_{v}$ for each $v\in
Q_{0}$ and define a linear function
$\Theta:\mathbb{Z}Q_{0}\rightarrow \mathbb{Z}$, by
$$\Theta(M):=\Theta(\underline{\dim}M)=\sum_{v\in Q_{0}}\Theta_{v}\dim_{k}M_{v}\; .$$

Let $\sigma$ be a set of strictly positive numbers $\sigma_{v}$
for each $v\in Q_{0}$, and define a (strictly positive) linear
function $\sigma:\mathbb{Z}Q_{0}\rightarrow \mathbb{Z}$, by
$$\sigma(M):=\sigma(\underline{\dim}M)=\sum_{v\in Q_{0}}\sigma_{v}\dim_{k}M_{v}\; .$$
We call $\sigma(M)$ \emph{the total dimension of $M$}. We will
refer to $\Theta$ and $\sigma$ indistinctly meaning the sets of
numbers or the linear functions.

For a non-zero representation $M$ of $Q$ over $k$, define its
slope by
$$\mu_{(\Theta,\sigma)}(M):=\frac{\Theta(M)}{\sigma(M)}\; .$$

\begin{dfn}
\label{Qstability}
A representation $M$ of $Q$ over $k$ is $(\Theta,\sigma)$-semistable if for all non-zero proper subrepresentations $M'\subset M$, we have
$$\mu_{(\Theta,\sigma)}(M')\leq\mu_{(\Theta,\sigma)}(M)\; .$$
If the inequality is strict for every non-zero proper subrepresentation, we say that $M$ is $(\Theta,\sigma)$-stable.
\end{dfn}

\begin{lem}
\label{changesthetasigma}
If we multiply the linear function $\Theta$ by a non-negative integer, or if we add an integer multiple of the strictly positive
linear function $\sigma$ to $\Theta$, the semistable (resp. stable) representations remain semistable (resp. stable). 
\end{lem}
\begin{pr}
Let $\Theta'=a\cdot \Theta+b\cdot \sigma,\; a,b\in\mathbb{Z},\; a>0$ be another linear function and note that
$$\frac{\Theta'(M')}{\sigma(M')}\leq \frac{\Theta'(M)}{\sigma(M)}\Leftrightarrow \frac{a\cdot \Theta(M')+b\cdot \sigma(M)}{\sigma(M')}\leq 
\frac{a\cdot \Theta(M)+b\cdot \sigma(M)}{\sigma(M)}$$
$$\Leftrightarrow \frac{\Theta(M')}{\sigma(M')}\leq \frac{\Theta(M)}{\sigma(M)}\; .$$\end{pr}

\begin{rem}
\label{alastair}
In \cite{Ki}, the stability condition (c.f
\cite[Definition 1.1]{Ki}) is formulated by not considering representations with different dimension vectors.
This leads to the construction of a moduli space and
$S$-filtrations (or Jordan-H\"{o}lder filtrations) but not to define a Harder-Narasimhan filtration,
for which is needed a slope condition as in Definition
\ref{Qstability}.

This slope stability condition, the $(\Theta,\sigma)$-stability
(c.f. Definition \ref{Qstability}), can be turned out into a
stability condition as in \cite{Ki}, by clearing denominators
$$\theta(M')=\Theta(M)\sigma(M') - \sigma(M)\Theta(M')\; ,$$
where $\theta$ is the function in \cite[Definition 1.1]{Ki} (observe that
$\theta(M)=0$), $\Theta$ and $\sigma$ are as in Definition
\ref{Qstability}, and $M'\subset M$ is a subrepresentation.

We will apply this in Proposition \ref{GITstab-stab}, to
relate $(\Theta,\sigma)$-stability with GIT stability.
\end{rem}

\begin{rem}
The definition of stability which appears in \cite{Re} considers
$\sigma_{v}=1$ for each $v\in Q_{0}$, although we consider a
strictly positive linear function $\sigma$ in general. The
notation of $\sigma$ agrees with \cite{AC}, \cite{ACGP},
\cite{Sch}, while $\Theta$ agrees with \cite{Re} but in the other
references it is substituted by different notations closer to
classical moduli problems where the stability notion depends on
parameters ($\tau$-stability or $\rho$-stability).
\end{rem}

\begin{lem}\cite[Definition 1]{Ru},
\cite[Lemma 4.1]{Re} \label{stabilitycategory} Let $0\rightarrow
X\rightarrow Y\rightarrow Z\rightarrow 0$ be a short exact
sequence of non-zero representations of $Q$ over $k$. Then
$\mu_{(\Theta,\sigma)}(X)< \mu_{(\Theta,\sigma)}(Y)$ if and only
if $\mu_{(\Theta,\sigma)}(X)<\mu_{(\Theta,\sigma)}(Z)$ if and only
if  $\mu_{(\Theta,\sigma)}(Y)<\mu_{(\Theta,\sigma)}(Z)$.
\end{lem}
\begin{pr}
 Note that $\sigma(Y)=\sigma(X)+\sigma(Z)$ and, therefore
$$\mu_{(\Theta,\sigma)}(Y)=\frac{\Theta(Y)}{\sigma(Y)}=\frac{\Theta(X)+\Theta(Z)}{\sigma(X)+\sigma(Z)}\; ,$$
from which the statement follows.
\end{pr}

\begin{thm}\cite[Theorem 2]{Ru},
\cite[Lemma 4.7]{Re}
\label{HN}
Given linear functions $\Theta$ and $\sigma$, (being $\sigma$ strictly positive), every representation $M$ of $Q$ over $k$ has a unique filtration
$$0\subset M_{1}\subset M_{2}\subset \ldots \subset M_{t}\subset M_{t+1}=M$$
verifying the following properties, where $M^{i}:=M_{i}/M_{i-1}$
\begin{enumerate}
\item $\mu_{(\Theta,\sigma)}(M^{1})>\mu_{(\Theta,\sigma)}(M^{2})>\ldots >\mu_{(\Theta,\sigma)}(M^{t})>\mu_{(\Theta,\sigma)}(M^{t+1})$
\item The quotients $M^{i}$ are $(\Theta,\sigma)$-semistable
\end{enumerate}
This filtration is called the \emph{Harder-Narasimhan filtration} of $M$ (with respect to $\Theta$ and $\sigma$).
\end{thm}
\begin{pr}
Using Lemma \ref{stabilitycategory} we can prove the existence of
a unique subrepresentation $M_{1}$, whose slope is maximal among
all the subrepresentations of $M$, and of maximal total dimension
$\sigma(M_{1})$ between those of maximal slope (c.f.
\cite[Proposition 1.9]{Ru}, \cite[Lemma 4.4]{Re}). Then, proceed
by recursion on the quotient $M/M_{1}$.
\end{pr}

\section{Moduli space of representations of quivers}
Fix $k$ an algebraically closed field of arbitrary characteristic.
Fix a dimension vector $d\in\mathbb{Z}Q_{0}$ and fix $k$-vector
spaces $M_{v}$ of dimension $d_{v}$ for all $v\in Q_{0}$. Fix
linear functions $\Theta, \sigma:\mathbb{Z}Q_{0}\rightarrow
\mathbb{Z}$, being $\sigma$ strictly positive. We recall the
construction by King (c.f. \cite{Ki}) of a moduli space for
representations of $Q$ over $k$ with dimension vector $d$.

Consider the affine $k$-space
$$\mathcal{R}_{d}(Q)=\bigoplus_{\alpha:v_{i}\rightarrow v_{j}}\Hom_{k}(M_{v_{i}},M_{v_{j}})$$
whose points parametrize representations of $Q$ on the $k$-vector
spaces $M_{v}$. The reductive linear algebraic group
$$G_{d}=\prod_{v\in Q_{0}}GL(M_{v})$$
acts on $\mathcal{R}_{d}(Q)$ by
$$(g_{v_{i}})_{v_{i}}\cdot (M_{\alpha})_{\alpha}=(g_{v_{j}}M_{\alpha}g_{v_{i}}^{-1})_{\alpha:v_{i}\rightarrow
v_{j}}$$ and the $G_{d}$-orbits of $M$ in $\mathcal{R}_{d}(Q)$
correspond bijectively to the isomorphism classes $[M]$ of
$k$-representations of $Q$ with dimension vector $d$. We will use
Geometric Invariant Theory to take the quotient of
$\mathcal{R}_{d}(Q)$ by $G_{d}$ and construct a moduli space of
representations of the quiver $Q$ on the $k$-vector spaces
$M_{v}$.

The action of $G_{d}$ on the affine space $\mathcal{R}_{d}(Q)$ can
be lifted by a character $\chi$ to the (necessarily trivial) line
bundle $L$ required by the Geometric Invariant Theory. Note that
the subgroup of the diagonal scalar matrices in $G_{d}$,
$$\Delta=\{(t1,\ldots,t1):t\in k^{\ast}\}\; ,$$ acts
trivially on $\mathcal{R}_{d}(Q)$. Then, we have to choose $\chi$
in such a way that $\Delta$ acts trivially on the fiber, in other
words, $\chi(\Delta)=1$.

Then, using the linear functions $\Theta$ and $\sigma$, consider
the character
$$\chi_{(\Theta,\sigma)}((g_{v})_{v}):=\prod_{v\in Q_{0}}\det(g_{v})^{(\Theta(d)\sigma_{v}-\sigma(d)\cdot \Theta_{v})}$$
of $G_{d}$, and note that $\chi_{(\Theta,\sigma)}(\Delta)=1$,
because $\sum_{v\in
Q_{0}}(\Theta(d)\sigma_{v}-\sigma(d)\Theta_{v})\cdot d_{v}=0$.

\begin{dfn}\cite[Definition 2.1]{Ki}
A point $x\in\mathcal{R}_{d}(Q)$ is $\chi$-semistable if there is
a relative invariant $f\in
k[\mathcal{R}_{d}(Q)]^{G_{d},\chi^{n}_{(\Theta,\sigma)}}$ with
$n\geq 1$, such that $f(x)\neq 0$.
\end{dfn}

The algebraic quotient will be given by
$$\mathcal{R}_{d}(Q)/\!\!/(G_{d},\chi_{(\Theta,\sigma)})=\Proj\big(\bigoplus_{n\geq
0}k[\mathcal{R}_{d}(Q)]^{G_{d},\chi^{n}_{(\Theta,\sigma)}}\big)\;
.$$

\begin{prop}
\label{GITstab-stab} A point $x_{M}\in \mathcal{R}_{d}(Q)$
corresponding to a representation $M\in \mod kQ$ is
$\chi_{(\Theta,\sigma)}$-semistable (resp.
$\chi_{(\Theta,\sigma)}$-stable) for the action of $G_{d}$ if and
only if $M$ is $(\Theta,\sigma)$-semistable (resp.
$(\Theta,\sigma)$-stable).
\end{prop}
\begin{pr}
It follows easily from \cite[Proposition 3.1]{Ki} and the
observation in Remark \ref{alastair}. In \cite{Ki}, given a linear
function $\theta$, a representation $M$ is $\theta$-semistable if
$\theta(M)=0$ and, for every subrepresentation $M'\subset M$, we
have $\theta(M')\geq 0$ (c.f. \cite[Definition 1.1]{Ki}). Then,
\cite[Proposition 3.1]{Ki} relates the $\theta$-stability with the
$\chi_{\theta}$-stability, where the character is
$$\chi_{\theta}((g_{v})_{v}):=\prod_{v\in Q_{0}}\det(g_{v})^{\theta_{v}}\; .$$
Hence, the $\chi_{(\Theta,\sigma)}$-stability with the character
given by
$$\chi_{(\Theta,\sigma)}((g_{v})_{v}):=\prod_{v\in Q_{0}}\det(g_{v})^{(\Theta(d)\sigma_{v}-\sigma(d)\Theta_{v})}\; ,$$
is equivalent to the $(\Theta,\sigma)$-stability in Definition
\ref{Qstability} because, given a subrepresentation $M'\subset M$,
the expression
$$\sum_{v\in
Q_{0}}(\Theta(M)\sigma_{v}-\sigma(M)\Theta_{v})\cdot \dim M'_{v}=\Theta(M)\sigma(M')-\sigma(M)\Theta(M')\geq 0$$
is equivalent to
$$\frac{\Theta(M')}{\sigma(M')}\leq \frac{\Theta(M)}{\sigma(M)}\; .$$\end{pr}

Now, denote by $\mathcal{R}^{(\Theta,\sigma)-ss}_{d}(Q)$ the set of
$\chi_{(\Theta,\sigma)}$-semistable points.

\begin{thm}\cite[Proposition 4.3]{Ki},
\cite[Corollary 3.7]{Re} The variety
$\mathfrak{M}^{(\Theta,\sigma)}_{d}(Q)=\mathcal{R}^{(\Theta,\sigma)-ss}_{d}(Q)/\!\!/(G_{d},\chi_{(\Theta,\sigma)})$
is a moduli space which parametrizes $S$-equivalence classes
of $(\Theta,\sigma)$-semistable representations of $Q$ of
dimension vector $d$. It is projective over the ordinary quotient $\mathcal{R}_{d}(Q)/\!\!/G_{d}$.
\end{thm}

If the quiver has no oriented cycles or it is chosen so that there is a unique semisimple representation for each dimension vector (e.g. the quiver associated to a finite
dimensional algebra over an algebraically closed field in \cite{Ki}), the ordinary quotient $\mathcal{R}_{d}(Q)/\!\!/G_{d}$ consists on one single point, 
hence the moduli space is projective. 

$$ $$

By the Hilbert-Mumford criterion we can characterize
$\chi_{(\Theta,\sigma)}$-semistable points by its behavior under
the action of $1$-parameter subgroups. A $1$-parameter subgroup of
$G_{d}=\prod_{v\in Q_{0}}GL(M_{v})$ is a non-trivial homomorphism
$\Gamma:k^{\ast}\rightarrow G_{d}$. There exist bases of the
vector spaces $M_{v}$ such that $\Gamma$ takes the diagonal form

$$\left(
  \begin{array}{ccc}
    t^{\Gamma_{v_{1},1}} &  & \\
     & \ddots & \\
     &  & t^{\Gamma_{v_{1},t_{1}+1}} \\
  \end{array}
\right)\times \cdots \times
\left(
  \begin{array}{ccc}
    t^{\Gamma_{v_{s},1}} &  & \\
     & \ddots & \\
     &  & t^{\Gamma_{v_{s},t_{s}+1}} \\
  \end{array}
\right)$$

where $v_{1},\ldots ,v_{s}\in Q_{0}$ are the vertices of the quiver.

Let $x\in\mathcal{R}_{d}(Q)$ and suppose that $\lim_{t\rightarrow
0}\Gamma\cdot x$ exists and is equal to $x_{0}$. Then $x_{0}$ is a
fixed point for the action of $\Gamma$, and $\Gamma$ acts on the
fiber of the trivial line bundle over $x_{0}$ as multiplication by
$t^{a}$. Define the following numerical function,
$$\mu_{\chi_{(\Theta, \sigma)}}(x,\Gamma)=-a\; .$$
The next proposition establishes the so-called ''Hilbert-Mumford
numerical criterion``:

\begin{prop}\cite[Proposition 2.5]{Ki}
\label{mrw} A point $x_{M}\in \mathcal{R}_{d}(Q)$ corresponding to a representation $M$ is
$\chi_{(\Theta,\sigma)}$-semistable if and only if every
$1$-parameter subgroup $\Gamma$ of $G_{d}$, for which
$\lim_{t\rightarrow 0} \Gamma(t)\cdot x_{M}$ exists, satisfies
$\mu_{\chi_{(\Theta, \sigma)}}(x_{M},\Gamma)\leq 0$.
\end{prop}

\begin{rem}
Note that in Proposition \ref{mrw} we change the sign of the
numerical function $\mu_{\chi_{(\Theta, \sigma)}}(x_{M},\Gamma)$
with respect to \cite{Ki} (as we did when changing the character in the proof of Proposition \ref{GITstab-stab}), in congruence with \cite{Ke} and
\cite{GSZ}.
\end{rem}

The action of a $1$-parameter subgroup $\Gamma$ of $G_{d}$
provides a decomposition of each vector space $M_{v}$ associated
to each vertex $v\in Q_{0}$, in weight spaces
$$M_{v}=\bigoplus_{n\in \mathbb{Z}}M_{v}^{n}$$
where $\Gamma(t)$ acts on the weight space $M_{v}^{n}$ as
multiplication by $t^{n}$. Every $1$- parameter subgroup, for
which $\lim_{t\rightarrow 0} \Gamma(t)\cdot x$ exists, determines
a weighted filtration $M_{\bullet}\subset M$ of subrepresentations
(c.f. \cite{Ki})
$$0\subset M_{1}\subset M_{2}\subset \ldots \subset M_{t}\subset M_{t+1}=M$$
where $M_{i}$ is the subrepresentation with vector spaces
$M_{v,i}:=\bigoplus _{n\leq i}M_{v}^{n}$ for each vertex $v\in
Q_{0}$, and the weight corresponding to each quotient
$M^{i}=M_{i}/M_{i-1}$ is $\Gamma_{i}$. Note that two $1$-parameter
subgroups giving the same filtration are conjugated by an element
of the parabolic subgroup of $G_{d}$ defined by the filtration.
Therefore, the numerical function $\mu_{\chi_{(\Theta,
\sigma)}}(x_{M},\Gamma)$, has a simple expression in terms of the
filtration $M_{\bullet}\subset M$ (c.f. calculation in \cite{Ki}):
\begin{equation}
\label{pairing} \mu_{\chi_{(\Theta, \sigma)}}(x_{M},\Gamma)=\sum_{v\in
Q_{0}}\big[\big(\Theta(M)\sigma_{v}-\sigma(M)\Theta_{v}\big)\cdot
\sum_{i=1}^{t_{v}+1}\Gamma_{v,i}\dim M_{v}^{i}\big]\; .
\end{equation}

Let $d_{i}$, $d^{i}$ be the dimension vectors of the
subrepresentation $M_{i}$ and the quotient $M^{i}=M_{i}/M_{i-1}$,
respectively. The action of $\Gamma$ on the point corresponding to
a representation $M$ has different weights for each vertex $v\in
Q_{0}$, but collect all different weights $\Gamma_{i}$ corresponding to any vertex and form
the vector
$$\Gamma=(\Gamma_{1},\Gamma_{2},\ldots,\Gamma_{t},\Gamma_{t+1})$$ verifying $\Gamma_{1}<\Gamma_{2}<\ldots<\Gamma_{t}<\Gamma_{t+1}$.
Hence, (\ref{pairing}) turns out to be
\begin{equation}
\label{pairing2} \mu_{\chi_{(\Theta,
\sigma)}}(x_{M},\Gamma)=\sum_{i=1}^{t+1}\Gamma_{i}\cdot
[\Theta(M)\cdot \sigma(M^{i})-\sigma(M)\cdot \Theta (M^{i})]\;,
\end{equation}
and Proposition \ref{mrw} can be rewritten in terms of filtrations of $M$.

\begin{prop}
\label{mrw2} A point $x_{M}\in \mathcal{R}_{d}(Q)$ corresponding
to a representation $M$ of $Q$ over $k$, is
$\chi_{(\Theta,\sigma)}$-semistable if and only if every
$1$-parameter subgroup $\Gamma$ of $G_{d}$, defining a filtration
of subrepresentations of $M$, $$0\subset M_{1}\subset M_{2}\subset
\ldots \subset M_{t}\subset M_{t+1}=M\; ,$$ satisfies that
$$\mu_{\chi_{(\Theta,
\sigma)}}(x_{M},\Gamma)=\sum_{i=1}^{t+1}\Gamma_{i}\cdot
[\Theta(M)\cdot \sigma(M^{i})-\sigma(M)\cdot \Theta (M^{i})]\leq
0\; .$$
\end{prop}

\section{Kempf theorem}

Given a weighted filtration of $M$,
$$0\subset M_{1}\subset M_{2}\subset \ldots \subset M_{t}\subset M_{t+1}=M$$
and $\Gamma_{1}<\Gamma_{2}<\ldots<\Gamma_{t}<\Gamma_{t+1}$, define the following function which we call the \emph{Kempf function},
\begin{equation}
\label{kempffunction}
K(x_{M},\Gamma)=\frac{\sum_{i=1}^{t+1}\Gamma_{i}\cdot
[\Theta(M)\cdot \sigma(M^{i})-\sigma(M)\cdot \Theta
(M^{i})]}{\sqrt{\sum_{i=1}^{t+1}\sigma(M^{i})\cdot
\Gamma_{i}^{2}}}
\end{equation}

We recall a theorem by Kempf (c.f. \cite[Theorem 2.2]{Ke}) stating
that whenever there exists any $\Gamma$ giving a positive value
for the numerator of the Kempf function, there exists a unique
parabolic subgroup containing a unique $1$-parameter subgroup in
each maximal torus, giving maximum in the Kempf function i.e.,
there exists a unique weighted filtration of $M$ for which the
Kempf function achieves its maximum.

The Kempf function (\ref{kempffunction}) which appears in
\cite[Theorem 2.2]{Ke} is a rational function whose numerator is
equal to the numerical function $\mu_{\chi_{(\Theta,
\sigma)}}(x_{M},\Gamma)$ and the denominator is the \emph{length} of
the $1$-parameter subgroup $\Gamma$. Given a reductive algebraic
linear group $G$, there is a notion of \emph{length} defined by Kempf
(c.f. \cite[pg. 305]{Ke}) in $\Gamma(G)$, the set of all
$1$-parameter subgroups.

\begin{dfn}
\label{length}
 A \emph{length} is a non-negative function $\|\;\|$
 on $\Gamma(G)$ with values on the real numbers, invariant by
conjugation by rational points of $G$, and such that for any
maximal torus $T\subset G$, there is a positive definite integral
valued form $(\cdot,\cdot)$ in $\Gamma(T)$ with
$(\Gamma,\Gamma)=\|\Gamma\|^{2}$, for any $\Gamma\in \Gamma(T)$.
\end{dfn}

If $G$ is simple, in characteristic zero all choices of length will be multiples of the Killing form in the
Lie algebra $\mathfrak{g}$ (note that in this case $\Gamma(G)\subset \mathfrak{g}$). For an almost simple
group in arbitrary characteristic (a group $G$ whose center $Z$ is finite and $G/Z$ is simple, e.g. $SL(N)$ in positive characteristic), all lengths differ also by a scalar.

However, in this case, the group is a product of general linear groups, which is not simple. Then, there are several simple factors in the group and we can take a 
different multiple of the Killing form for each factor. Hence, observe that in the Kempf function (\ref{kempffunction}), the denominator of the expression is a
function verifying the properties of the definition of a length (c.f. Definition \ref{length}). The different multiples we take for each factor are the integer coefficients
of the strictly positive linear function $\sigma$. 

Therefore, we can rewrite \cite[Theorem 2.2]{Ke} in our case as follows:

\begin{thm}
\label{Kempf} Given a $\chi_{(\Theta,\sigma)}$-unstable point
$x_{M}\in \mathcal{R}_{d}(Q)$ corresponding to a representation $M$, there exists a unique weighted
filtration, i.e. $0\subset M_{1}\subset \cdots \subset M_{t+1}= M$
and real numbers
$\Gamma_{1}<\Gamma_{2}<\ldots<\Gamma_{t}<\Gamma_{t+1}$, called the
\emph{Kempf filtration of M}, such that the Kempf function $K$
achieves the maximum among all filtrations and weights verifying
$\Gamma_{1}<\Gamma_{2}<\ldots<\Gamma_{t}<\Gamma_{t+1}$.
\end{thm}

Note that the length we are considering depends on the choice of
$\sigma$ and the Kempf function depends both on the length and the
linearization of the group action, hence depends both on $\Theta$
and $\sigma$. In order to relate the Kempf filtration
of $M$ with the Harder-Narasimhan filtration, which also depends on
$\Theta$ and $\sigma$, we dispose the parameters
conveniently.

\section{Results on convexity}

Next, we prove a result about convexity for functions which are similar to the Kempf function. The vector which maximizes such functions verifies some properties that will be strongly
related to the properties of the Harder-Narasimhan filtration. In this section we recall the results of \cite[Section 2]{GSZ}.

Consider $\mathbb{R}^{t+1}$ together with an inner product
$(\cdot,\cdot)$ defined by the  diagonal matrix
 $$
 \left(
 \begin{array}{ccc}
 b^1 & & 0 \\
  & \ddots & \\
 0 & & b^{t+1}\\
 \end{array}
 \right)
 $$
where $b^i$ are positive integers. Let
$$\mathcal{C}= \big\{ x\in \mathbb{R}^{t+1} : x_1<x_2<\cdots <x_{t+1} \big\}\; ,$$
and $v= (v_1,\cdots,v_{t+1})\in \mathbb{R}^{t+1}$ verifying $\sum_{i=1}^{t+1}b^{i}v_{i}=0$. Define the function
\begin{eqnarray*}
  \label{eq:mu}
\mu_{v}:\overline{\mathcal{C}} & \to & \mathbb{R}\\
\Gamma & \mapsto & \mu_{v}(\Gamma)=\frac{(\Gamma,v)}{||\Gamma||}
\end{eqnarray*}
and note that  $\mu_{v}(\Gamma)=||v||\cdot \cos(\Gamma,v)$.

We assume that there exists $\Gamma\in\overline{\mathcal{C}}$ with
$\mu_v(\Gamma)>0$ and then, we would like to find a vector
$\Gamma\in \overline{\mathcal{C}}$ which maximizes the function
$\mu_{v}$. Define $w^{i}=-b^{i}\cdot v_{i}$,
$w_{i}=w^{1}+\cdots+w^{i}$, $b_{i}=b^{1}+\cdots+b^{i}$ and draw a
graph joining the points with coordinates $(b_{i},w_{i})$, each
segment having slope $-v_{i}$. Now draw the convex envelope of
this graph (thick line in Figure $1$), denoting its coordinates by
$(b_{i},\widetilde{w_{i}})$, and define
$$\Gamma_{i}=-\frac{\widetilde{w_{i}}-\widetilde{w_{i-1}}}{b^{i}}\; .$$
In other words, the vector
$\Gamma_{v}=(\Gamma_1,\cdots,\Gamma_{t+1})$ defined in this way,
verifies that the quantities $-\Gamma_i$ are the slopes of the
convex envelope graph defined by $v$.

\setlength{\unitlength}{1cm}
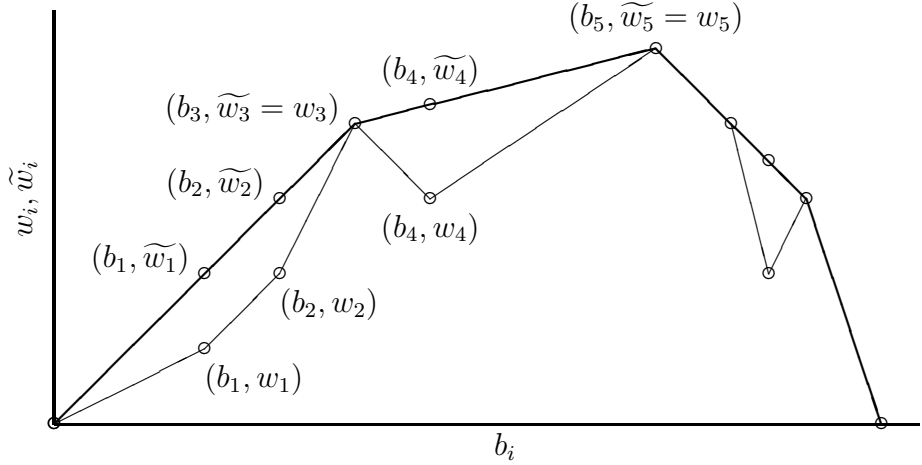
\begin{figure}[h]
\begin{picture}(13,7)(-1,-1)
\thicklines
\put(0,0){\line(1,0){11.5}}
\put(0,0){\line(0,1){5.5}}
\put(6,-0.3){\makebox(0,0)[c]{$b_i$}}
\put(-0.4,3){\makebox(0,0)[c]{\rotatebox{90}{$w_i,\widetilde{w_i}$}}}
\put(0,0){\makebox(0,0){$\circ$}}
\put(0,0){\line(1,1){4}}
\put(2,2){\makebox(0,0){$\circ$}}
\put(3,3){\makebox(0,0){$\circ$}}
\put(4,4){\makebox(0,0){$\circ$}}
\put(4,4){\line(4,1){4}}
\put(5,4.25){\makebox(0,0){$\circ$}}
\put(8,5){\makebox(0,0){$\circ$}}
\put(8,5){\line(1,-1){2}}
\put(9.5,3.5){\makebox(0,0){$\circ$}}
\put(10,3){\makebox(0,0){$\circ$}}
\put(10,3){\line(1,-3){1}}
\put(11,0){\makebox(0,0){$\circ$}}
\thinlines
\put(0,0){\makebox(0,0){$\circ$}}
\put(0,0){\line(2,1){2}}
\put(2,1){\makebox(0,0){$\circ$}}
\put(2,1){\line(1,1){1}}
\put(3,2){\makebox(0,0){$\circ$}}
\put(3,2){\line(1,2){1}}
\put(4,4){\line(1,-1){1}}
\put(5,3){\makebox(0,0){$\circ$}}
\put(5,3){\line(3,2){3}}
\put(9,4){\makebox(0,0){$\circ$}}
\put(9,4){\line(1,-4){0.5}}
\put(9.5,2){\makebox(0,0){$\circ$}}
\put(9.5,2){\line(1,2){0.5}}
\put(2,0.6){\makebox(0,0)[l]{$(b_1,w_1)$}}
\put(1.8,2.2){\makebox(0,0)[r]{$(b_1,\widetilde{w_1})$}}
\put(3,1.6){\makebox(0,0)[l]{$(b_2,w_2)$}}
\put(2.8,3.2){\makebox(0,0)[r]{$(b_2,\widetilde{w_2})$}}
\put(3.8,4.2){\makebox(0,0)[r]{$(b_3,\widetilde{w_3}=w_3)$}}
\put(5,4.7){\makebox(0,0)[c]{$(b_4,\widetilde{w_4})$}}
\put(5,2.6){\makebox(0,0)[c]{$(b_4,w_4)$}}
\put(8,5.4){\makebox(0,0)[c]{$(b_5,\widetilde{w_5}=w_5)$}}
\end{picture}
\caption{Convex envelope $\Gamma_{v}$ of the graph defined by $v$}
\end{figure}

\begin{thm}\cite[Theorem 2.2]{GSZ}
\label{max}
The vector $\Gamma_{v}$ defined in this way gives the maximum for the function $\mu_v$ on its domain.
\end{thm}

\section{Kempf filtration is Harder-Narasimhan filtration }

Finally, we study the geometrical properties of the Kempf filtration by associating to it a graph which encodes the two properties satisfied by
the Harder-Narasimhan filtration.

Let $\Theta:\mathbb{Z}Q_{0}\rightarrow \mathbb{Z}$ be a linear
function and let $\sigma:\mathbb{Z}Q_{0}\rightarrow \mathbb{Z}$ be
a strictly positive linear function. Let $M$ be a representation
of $Q$ over an algebraically closed field $k$ of arbitrary
characteristic, which is  $(\Theta,\sigma)$-unstable. Consider the
$\chi_{(\Theta,\sigma)}$-unstable point $x_{M}\in
\mathcal{R}_{d}(Q)$ associated to $M$, by Proposition
\ref{GITstab-stab}. Let $0\subset M_{1}\subset \cdots \subset
M_{t+1}= M$ and
$\Gamma_{1}<\Gamma_{2}<\ldots<\Gamma_{t}<\Gamma_{t+1}$ be the
Kempf filtration of $M$, by Theorem \ref{Kempf}.

Let $M^{i}=M_{i}/M_{i-1}$ be the quotients of the filtration.
Consider the inner product in $\mathbb{R}^{t+1}$ given by the
matrix $$
 \left(
 \begin{array}{ccc}
 \sigma(M^1) & & 0 \\
  & \ddots & \\
 0 & & \sigma(M^{t+1})\\
 \end{array}
 \right)
 $$
where $\sigma(M^{i})>0$.
\begin{dfn}
\label{graph} Given a filtration $0\subset M_{1}\subset \cdots
\subset M_{t+1}= M$ of subrepresentations of $M$, define
$v=(v_{1},...,v_{t+1})$, where
$v_{i}=\Theta(M)-\frac{\sigma(M)}{\sigma(M^{i})}\Theta(M^{i})$,
the vector associated to the filtration.
\end{dfn}

Now we can identify the Kempf function (\ref{kempffunction}) with the function
in Theorem \ref{max},
$$K(x_{M},\Gamma)=\frac{\sum_{i=1}^{t+1}\Gamma_{i}\cdot [\Theta(M)\sigma(M^{i})-\sigma(M)\Theta(M^{i})]}{\sqrt{\sum_{i=1}^{t+1}\sigma(M^{i})\cdot
\Gamma_{i}^{2}}}=$$
$$=\frac{\sum_{i=1}^{t+1}\sigma(M^{i})\Gamma_{i}\cdot [\Theta(M)-\frac{\sigma(M)}{\sigma(M^{i})}\Theta (M^{i})]}{\sqrt{\sum_{i=1}^{t+1}
\sigma(M^{i})\cdot
\Gamma_{i}^{2}}}=\frac{(\Gamma,v)}{\|\Gamma\|}=\mu_{v}(\Gamma)\;
.$$ Note that $\sum_{i=1}^{t+1}b^{i}v_{i}=0$.

\begin{lem}\cite[Lemma 3.4, Lemma 3.5]{GSZ}
\label{lemmaAB}
Let $0\subset M_{1}\subset \cdots \subset M_{t+1}=M$ be the Kempf filtration of $M$ (cf. Theorem \ref{Kempf}).
Let $v=(v_{1},...,v_{t+1})$ the vector associated to this filtration by Definition \ref{graph}. Then,
\begin{enumerate}
\item The coordinates of $v$ verify $v_{1}<v_{2}<\ldots<v_{t}<v_{t+1}$
i.e., the graph of $v$ is convex.
\item The vector $v$ is the convex envelope of every refinement.
\end{enumerate}
\end{lem}

\begin{thm}
\label{KempfHN}
The Kempf filtration of $M$ is the Harder-Narasimhan filtration of $M$.
\end{thm}
\begin{pr}
The vector $v$ associated to the Kempf filtration of $M$ verifies properties $(1)$ and $(2)$ in
Lemma \ref{lemmaAB}, which are precisely the properties
$(1)$ and $(2)$ in Theorem \ref{HN}, respectively. By uniqueness
of the Harder-Narasimhan filtration of $M$, both filtrations do
coincide.
\end{pr}


\begin{thebibliography}{EMG}

\bibitem[AC]{AC}{L. \'Alvarez-C\'onsul, }
\textit{Some results on the moduli spaces of quiver bundles, } Geom. Dedicata
\textbf{139} (2009), 99-120.

\bibitem[ACGP]{ACGP}{L. \'Alvarez-C\'onsul, O. Garc\'ia-Prada, }
\textit{Hitchin-Kobayashi correspondence, quivers, and vortices, }
Comm. Math. Phys. \textbf{238} (2003), 1-33.

\bibitem[GSZ]{GSZ}{T. G\'omez, I. Sols, A. Zamora, }
\textit{A GIT characterization of the Harder-Narasimhan filtration, }
arxiv:1112.1886v2, (Preprint 2011).

\bibitem[He]{He}{W. H. Hesselink, }
\textit{Uniform instability in reductive groups, } J. Reine Angew.
Math. \textbf{304} (1978), 74-96.

\bibitem[Ho]{Ho}{V. Hoskins, }
\textit{Stratifications associated to reductive group actions on affine spaces, }
arxiv:1210.6811, (Preprint 2012).

\bibitem[Ke]{Ke}{G. Kempf, }
\textit{Instability in invariant theory, } Ann. Math. (2)
\textbf{108} no. 1 (1978), 299-316.

\bibitem[Ki]{Ki}{A. King, }
\textit{Moduli of representations of finite dimensional algebras, }
Q. J. Math. Oxf. Ser. (2), \textbf{45} (1994), 180, 515-530.

\bibitem[Re]{Re}{M. Reineke, }
\textit{Moduli of representations of quivers, }
arxiv:0802.2147v1, (Preprint 2008)

\bibitem[Ru]{Ru}{A. Rudakov, }
\textit{Stability for an abelian category, }
J. Algebra, \textbf{197} (1997), 231-245.

\bibitem[Sch]{Sch}{A.H.W. Schmitt, }
\textit{Moduli for decorated tuples of sheaves and representation spaces for quivers, }
Proc. Indian Acad. Sci. Math. Sci. \textbf{115} (2005) 15-49.

\bibitem[Tu]{Tu}{L. Tur, }
\textit{Dualite etrange sur le plan projectif, Ph.D. Thesis, Universit\'{e} de Nice-Sophia Antipolis, }
(2003).

\end{thebibliography}
\end{document}